\documentclass{amsart}

\usepackage{amsthm}
\usepackage{amsfonts}
\usepackage{amsmath}
\usepackage{amssymb}
\usepackage{mathrsfs}
\usepackage{fullpage}
\usepackage{stmaryrd}
\usepackage{hyperref}
\usepackage{verbatim}

\theoremstyle{plain}

\theoremstyle{definition}

\theoremstyle{remark}

\newcommand{\Begin}[2]{\begin{#1}\label{#2}}

\newcommand{\bSigma}{\mathbf{\Sigma}}

\newcommand{\bbC}{\mathbb{C}}

\newcommand{\bbP}{\mathbb{P}}
\newcommand{\bbQ}{\mathbb{Q}}
\newcommand{\bbR}{\mathbb{R}}

\newcommand{\forces}{\Vdash}
\newcommand{\analytic}{{\bSigma_1^1}}

\newcommand{\ZFC}{\mathsf{ZFC}}
\newcommand{\ZF}{\mathsf{ZF}}

\newcommand{\reals}{\bbR}

\newcommand{\OD}{\mathrm{OD}}

\newcommand{\AD}{\mathsf{AD}}
\newcommand{\AC}{\mathsf{AC}}
\newcommand{\DC}{\mathsf{DC}}
\newcommand{\ON}{\mathrm{ON}}

\begin{document}

\title{The Destruction of the Axiom of Determinacy by Forcings on $\reals$ When $\Theta$ Is Regular}

\author{William Chan}
\address{Department of Mathematics, University of North Texas, Denton, TX 76203}
\email{William.Chan@unt.edu}

\author{Stephen Jackson}
\address{Department of Mathematics, University of North Texas, Denton, TX 76203}
\email{Stephen.Jackson@unt.edu}

\begin{abstract}
$\ZF + \AD$ proves that for all nontrivial forcings $\bbP$ on a wellorderable set of cardinality less than $\Theta$, $1_\bbP \forces_\bbP \neg\AD$. $\ZF + \AD + \Theta$ is regular proves that for all nontrivial forcing $\bbP$ which is a surjective image of $\reals$, $1_\bbP \forces_\bbP \neg\AD$. In particular, $\mathsf{ZF + \AD + V = L(\reals)}$ proves that for every nontrivial forcing $\bbP \in L_\Theta(\reals)$, $1_\bbP \forces_\bbP \neg\AD$.
\end{abstract}

\thanks{March 16, 2019. The first author was supported by NSF grant DMS-1703708. The second author was supported by NSF grant DMS-1800323.}

\maketitle


\section{Introduction}\label{introduction}

Paul Cohen \cite{Independence-Continuum-Hypothesis} developed forcing which is a very flexible method of extending models of certain axioms of set theory (such as $\ZF$ or $\ZFC$) so that the resulting structures continue to satisfy these axioms of set theory. This technique has become a powerful tool for showing statements are independent of $\ZFC$. For example by \cite{Independence-Continuum-Hypothesis}, if $\ZFC$ is consistent, then $\ZFC + 2^{\aleph_0} = \aleph_1$ and $\ZFC + 2^{\aleph_0} > \aleph_1$ are both consistent. 

Descriptive set theory is the study of the definable aspects of mathematics. Various interesting properties are commonly considered while employing definable techniques to study mathematical problems. Some of these include the perfect set property, Lebesgue measurability, the property of Baire, partition relations on ordinals, and certain properties of Turing degrees. These properties in their full generality are all incompatible with $\AC$, the axiom of choice.

These properties are interesting and appeared naturally in classical descriptive set theory. Definable instances of these properties are provable in more basic axiom systems such as $\ZF$, in the same way that definable instances of the axiom of choice, for example, coanalytic uniformization, is provable in $\ZF$. This analogy justifies the study of the consequences of the full generalization of these properties just as one does with $\AC$, the full generalization of definable selection principles. 

The axiom of determinacy, $\AD$, has developed into a comprehensive framework for studying the properties mentioned above in their full generality. As customary in descriptive set theory, $\reals$ will denote the Baire space, ${}^\omega\omega$, of functions from $\omega$ into $\omega$. For each $A \subseteq \reals$, let $G_A$ be the game where Player 1 and 2 take turns playing elements of $\omega$. After infinitely many stages, a single $f \in \reals$ has been produced. Player 1 wins this run of $G_A$ if and only if $f \in A$. The axiom of determinacy states that for all $A \subseteq \reals$, one of the two players has a winning strategy for $G_A$. $\AD$ implies the perfect set property, Lebesgue measurability, Baire property for all sets of reals, and there are many cardinals with various partition properties. As with all these properties, definable fragments of $\AD$ can be proved in $\ZF$, for example, Martin showed all games $G_A$ where $A$ is Borel are determined under $\ZF$.

One can wonder if the forcing construction which has been fruitful for studying consistency results over $\ZFC$ could be useful for $\AD$. The most basic question would be to understand what forcings over $\AD$ could preserve $\AD$. By the nature of $\AD$, if one does not change $\reals$ or $\mathscr{P}(\reals)$, then $\AD$ will be preserved. Therefore the question becomes what forcings which disturb $\reals$ or $\mathscr{P}(\reals)$ can still preserve $\AD$.

Ikegami and Trang initiated the study of the preservation of $\AD$ under forcing. They showed that many forcings, such as Cohen forcing, can never preserve $\AD$. They also showed that if one is working with natural models of $\AD$, i.e. models satisfying $\mathsf{ZF + AD^+ + V = L(\mathscr{P}(\reals))}$, then any forcing which preserve $\AD$ must preserve $\Theta$, where $\Theta$ is the supremum of the ordinals which are surjective images of $\reals$. They also showed that the consistency of $\ZF + \AD^+ + \Theta > \Theta_0$ implies the consistency of $\ZF + \AD$ and there is a forcing which preserve $\AD$ and increases $\Theta$. Thus necessarily this forcing must disturb $\mathscr{P}(\reals)$ by adding a new set of reals. 

The following are some examples of concrete forcings applied within $\AD$. They all destroy the axiom of determinacy for various reasons. These examples give some emperical evidence that most small forcings can not preserve $\AD$ and also motivate the general arguments presented throughout the paper.

Let $\bbC$ denote Cohen forcing. Cohen forcing adds a generic filter which is equiconstrucible from a generic real, called the Cohen generic real. Ikegami and Trang observed that if $g$ is a Cohen generic real over $V$, then $V[g] \models$ ``$\reals^V$ does not have the Baire property''. Hence $V[g] \models \neg\AD$. 

Note that Woodin has shown that if $V \models \ZFC$ has a proper class of Woodin cardinals, then for any $\bbP$ and $G \subseteq \bbP$ which is $\bbP$-generic over $V$, $L(\reals)^V$ is elementarily equivalent to $L(\reals)^{V[G]}$. This setting implies that $L(\reals)^V \models \AD$. Let $g$ be a Cohen real which is generic over $V$. Note $g \in L(\reals)^{V[g]}$ and by Woodin's result, $L(\reals)^{V[g]} \models \AD$. However, $L(\reals)[g] \models \neg \AD$ by the result of Ikegami and Trang of the previous paragraph. Observe that the elements of the ground model always belong to its forcing extension. Thus $\reals^V \in L(\reals)[g]$; however, $\reals^V \notin L(\reals)^{V[g]}$.

Assume $\ZF + \DC_\reals + \AD$. Let $\mathrm{Coll}(\omega_1,\omega_2)$ be the forcing consisting of countable partial functions from $\omega_1$ into $\omega_2$ ordered by reverse extension. By $\AD$ and the coding lemma, there is a surjection $\pi : \reals \rightarrow \mathrm{Coll}(\omega_1,\omega_2)$. Suppose there was a $G \subseteq \mathrm{Coll}(\omega_1,\omega_2)$ generic over $V$ such that $V[G] \models \AD$. Since $\mathrm{Coll}(\omega_1,\omega_2)$ is countably closed and $\DC_\reals$ holds, no new reals are added. $V[G]$ has a surjection of $\omega_1^V$ onto $\omega_2^V$. Thus there is a new subset of $\omega_1^V$ which codes an ordering of $\omega_1^V$ of length $\omega_2^V$. In $V$, let $\pi : \reals^V \rightarrow \omega_1^V$ be a surjection. By the coding lemma in $V[G] \models \AD$, there is some real which codes this new subset of $\omega_1^V$ with respect to $\pi$. This is impossible if there are no new reals. Thus $V[G]$ can not satisfy $\AD$. 

Suppose $\kappa$ is a cardinal. The partition relation $\kappa \rightarrow (\kappa)^\lambda_2$ is the statement that for all $\Phi : [\kappa]^\lambda \rightarrow 2$, there is a club $C \subseteq \kappa$ and an $i \in 2$ so that $\Phi(f) = i$ for all $f \in [\kappa]^\lambda$ of the correct type. The notion of correct type will be defined below and is needed to obtain a club set which is homogeneous. Martin showed that $\omega_1 \rightarrow (\omega_1)^{\omega_1}_2$ holds under $\AD$. 

Assume $\ZF + \AD$. Let $\mathrm{Coll}(\omega,\omega_1)$ be the forcing consisting of finite partial functions from $\omega$ into $\omega_1$ ordered by reverse extension. Suppose there is a $G \subseteq \mathrm{Coll}(\omega,\omega_1)$ generic over $V$ such that $V[G] \models \AD$. One can show that $\omega_1^{V[G]} = \omega_2^V$. Since $|\mathrm{Coll}(\omega,\omega_1)|^V = \aleph_1^V$, one can show that for all club $D \subseteq \omega_1^{V[G]} = \omega_2^V$, there is a club $C \subseteq \omega_2^V$ which belong to $V$ so that $V[G] \models C \subseteq D$. (This is shown later as the ground club property.) In $V$, let $\Phi : [\omega_2]^{\omega_2} \rightarrow 2$ be an arbitrary partition. Since $\Phi \in V[G]$ and $([\omega_2]^{\omega_2})^V \in V[G]$, within $V[G]$, define $\tilde \Phi : [\omega_1^{V[G]}]^{\omega_1^{V[G]}} \rightarrow 2$ by
$$\tilde\Phi(f) = 
\begin{cases}
0 & \quad f \in ([\omega_2]^{\omega_2})^V \wedge \Phi(f) = 0 \\
1 & \quad \text{otherwise}
\end{cases}.$$
If $V[G] \models \AD$, then $\omega_1^{V[G]} \rightarrow (\omega_1^{V[G]})^{\omega_1^{V[G]}}$ implies there is some club $D \subseteq \omega_2^V$ so that $D$ is homogeneous for $\tilde \Phi$. Let $C \in V$ be club in $\omega_2^V$ so that $V[G] \models C \subseteq D$. One can show that $C$ is homogenous for $\Phi$ in $V$. Since $\Phi$ was an arbitrary partition, one has established $\omega_2 \rightarrow (\omega_2)^{\omega_2}_2$ in $V$. But Martin and Paris showed that $\AD$ implies $\omega_2 \rightarrow (\omega_2)^{\omega_2}_2$ is not true. Contradiction. So $\mathrm{Coll}(\omega,\omega_1)$ can not preserve $\AD$.

Revisiting the Cohen forcing $\bbC$: Assume $\ZF + \AD$. Suppose there was a Cohen generic real $g$ over $V$ such that $V[g] \models \AD$. Since $|\bbC| = \aleph_0$, every $D \subseteq \omega_1$ in $V[g]$ has a $C \in V$ which is a club subset of $\omega_1$ so that $V[g] \models C \subseteq D$. (Again this is the ground club property.) Note that $([\omega_1]^{\omega_1})^V \in V[g]$, so one may define a function $\Phi : [\omega_1]^{\omega_1} \rightarrow 2$ in $V[g]$ as follows:
$$\Phi(f) = \begin{cases}
0 & \quad f \in ([\omega_1]^{\omega})^V \\
1 & \quad \text{otherwise}
\end{cases}.$$
$V[g] \models \AD$, so by $\omega_1 \rightarrow (\omega_1)^\omega_2$, let $D \subseteq \omega_1$ be a club set homogeneous for $\Phi$. Let $C \subseteq D$ be a club in $V$ so that $V[g] \models C \subseteq D$. Taking any $f \in ([C]^{\omega_1})^V$ of the correct type, one can show that in $V[g]$, $C$ is homogeneous for $\Phi$ taking value $0$. Let $c_i$ denote the $(\omega \cdot i + \omega)^\text{th}$ element of $C$. As $C \in V$, $\langle c_i : i \in \omega \rangle \in V$. Pick $z \in \reals^{V[g]}$. Let $f_z \in ([C]^\omega)^{V[g]}$ be defined by letting $f_z$ be the increasing enumeration of $\{c_i : i \in z\}$. The function $f_z$ is of the correct type so $\Phi(f_z) = 0$. Thus $f_z \in V$. Since $z = \{i \in \omega : c_i \in f_z\}$, one has that $z \in V$. It has been shown that $\reals^V = \reals^{V[g]}$ which is impossible since $g \in \reals^{V[g]} \setminus \reals^V$. 

These examples suggest that ``small'' nontrivial forcings should not be able to preserve $\AD$. The examples also seem to indicate that the partition property and the ground model club phenomenon appears to be common aspects of these arguments. 

The axiom of determinacy by its definition influences the sets which are surjective images of $\reals$. It is reasonable to ask whether a nontrivial forcing which itself is within the realm of determinacy (i.e. is a surjective image of $\reals$) must disturb $\reals$ or $\mathscr{P}(\reals)$ and if so, can it preserve $\AD$. More specifically, if $V \models \AD$, $L(\reals)$ is the smallest model of determinacy containing $\reals^V$. One can ask if in $L(\reals)$, which is the most natural model of $\AD$, can a nontrivial forcing within the realm of determinacy, i.e. in $L_\Theta(\reals)$, preserve $\AD$. The following are the main questions:

\Begin{question}{main questions}
Assume $\ZF + \AD$. If $\bbP$ is a nontrivial forcing which is a surjective image of $\reals$, is it possible that $1_\bbP \forces_\bbP \AD$? 

Assume $\ZF + \AD + V = L(\reals)$. Is there any nontrivial $\bbP$ which is a surjective image of $\reals$ so that $1_\bbP \forces_\bbP \AD$?
\end{question}

The first question will be answered negatively if the assumptions are augmented with the condition that $\Theta$ is regular. Since $\Theta$ is regular in $L(\reals)$, this immediately gives the negative answer to the second question. The results of the paper are the following:
\\*
\\*\textbf{Theorem \ref{WO forcing below theta destroy AD}.}
\textit{Assume $\mathsf{ZF + AD}$. If $\bbP$ is a nontrivial wellorderable forcing of cardinality less than $\Theta$, then $1_{\bbP} \forces_{\bbP} \neg \AD$.}
\\*
\\*\indent The argument of the above theorem serves as a template for the main result. Its proof is a generalization of the example involving Cohen forcing. In discussion with Goldberg, a stronger result for wellorderable forcing can be shown using different techniques:
\\*
\\*\textbf{Corollary \ref{WO forcing adding new real destroy AD}.}
\textit{Assume $\ZF + \AD$. If $\bbP$ is a wellorderable forcing which adds a new real, the $1_\bbP \forces_\bbP \neg \AD$.}
\\*
\\* The main results are:
\\*
\\*\textbf{Theorem \ref{destruction AD when theta regular}.}
\textit{Assume $\ZF + \AD +$ $\Theta$ is regular. Suppose $\bbP$ is a nontrivial forcing which is a surjective image of $\reals$. Then $1_\bbP \forces_\bbP \neg\AD$.}
\\*
\\*\textbf{Corollary \ref{destruction AD natural model AD+}.}
\textit{Assume $\mathsf{ZF + AD + V = L(\reals)}$. No nontrivial forcing $\bbP \in L_\Theta(\reals)$ can preserve $\AD$.}

\textit{In fact, assume $\mathsf{ZF + AD^+ + \neg\AD_\reals + V = L(\mathscr{P}(\reals))}$. No nontrivial forcing which is the surjective image of $\reals$ can preserve $\AD$.}

\section{Ground Club Property}\label{ground club property}
Recall that if $A \subseteq \reals \times \reals^n$ and $e \in \reals$, $A_e = \{x \in \reals^n : (e,x) \in A\}$. 

\Begin{fact}{moschovakis coding lemma}
(Moschovakis) Assume $\ZF + \AD$. Let $\Gamma$ be a nonselfdual pointclass closed under continuous substitution, $\exists^\reals$, $\wedge$, and $\analytic \subseteq \Gamma$. Let $\prec \in \Gamma$ be a strict prewellordering. For each $a \in \mathrm{dom}(\prec)$, let $Q_a = \{b \in \mathrm{dom}(\prec) : a \preceq b \wedge b \preceq a\}$. Let $U \subseteq \reals^3$ be a $\Gamma$-universal set for subsets of $\reals^2$ in $\Gamma$. Let $Z \subseteq \mathrm{dom}(\prec) \times \reals$. Then there is an $e \in \reals$ so that 

\noindent (1) $U_e \subseteq Z$.

\noindent (2) For all $a \in \mathrm{dom}(\prec)$, $(U_{e})_{a} \neq \emptyset$ if and only if $Z_a \neq \emptyset$. 
\end{fact}

\begin{proof}
See \cite{Descriptive-Set-Theory} Section 7D.
\end{proof}

\Begin{fact}{coding sets of ordinals}
Assume $\ZF + \AD$. Let $X \subseteq \reals$ and $\pi : X \rightarrow \kappa$ be a surjection. Let $\prec$ be a strict prewellordering on $X$ defined by $x \prec y$ if and only if $\pi(x) < \pi(y)$. Let $\Gamma$ be a nonselfdual pointclass closed under continuous substitution, $\exists^\reals$, $\wedge$, and $\analytic \subseteq \Gamma$. Let $U$ be a fixed $\Gamma$-universal set for subsets of $\reals^2$ in $\Gamma$. For each $e \in \reals$, let $S_e^\pi = \{\alpha < \kappa : (\exists a)(\pi(a) = \alpha \wedge U_e(a,0)\}$. 

For all $C \subseteq \kappa$, there is some $e \in \reals$ so that $S_e^\pi = C$. 
\end{fact}

\begin{proof}
Let $Z = \{(a,0) : a \in X \wedge \pi(a) \in C\}$. Apply Fact \ref{moschovakis coding lemma}.
\end{proof}

\Begin{definition}{stable ordinals}
Assume $\ZF + \AD$. Let $A \subseteq \reals$. Let $\delta_A$ be the least ordinal $\delta$ so that $L_\delta(A,\reals) \prec_1 L(A,\reals)$, where $\prec_1$ denotes $\Sigma_1$ elementarity in a language that includes a predicate $\dot A$ and $\dot \reals$, which are always interpreted as $A$ and $\reals$, respectively. It is also the least ordinal $\delta$ so that $L_\delta(A,\reals)$ is an elementary substructure of $L(A,\reals)$ with respect $\Sigma_1$ formulas in the above language using elements of $\reals$, $\reals$ itself, and $A$ as parameters.

Let $\Sigma_1(L(A,\reals),\reals \cup \{\reals, A\})$ be the collection of sets in $L(A,\reals)$ which are $\Sigma_1$ definable in $L(A,\reals)$ using elements of $\reals$, $\reals$ itself, and $A$ as parameters.
\end{definition}

\Begin{definition}{prewellordering of stable}
Following \cite{Large-Cardinal-from-Determinacy} Section 2.4 and 2.5, the following is an explicit prewellordering of a subset of $\reals$ of length $\delta_A$ which is $\Sigma_1(L(A,\reals),\reals \cup \{\reals,A\})$:  

Let $T$ be the theory consisting of $\ZF$ without the power set axiom, ``$\reals$ exists'', and countable choice for $\reals$.

Let $\varphi_A(x,A,\dot\reals)$ denote a $\Sigma_1$ formula that defines the $\Sigma_1(L(A,\dot\reals), \{A,\dot\reals\})$ set, denoted $U_A$, which is universal for $\Sigma_1(L(A,\dot\reals), \dot\reals \cup \{A,\dot \reals\})$. For $x \in U_A$, let $\Theta_x$ be the least ordinal so that $L_{\Theta_x}(A,\dot\reals) \models T$ and $L_{\Theta_x}(A,\dot\reals) \models \varphi_A(x,A,\dot\reals)$. Define $\tilde\rho_A(x) = (\delta_A)^{L_{\Theta_x}(A,\dot\reals)}$. Let $\iota_A : \tilde\rho_A[U_A] \rightarrow \delta_A$ be the transitive collapse of $\tilde\rho[U_A]$. Let $\rho_A = \iota_A \circ \tilde\rho_A$. $\rho_A$ is a $\Sigma_1(L(A,\dot\reals),\{A,\dot\reals\})$ surjection of $U_A$ onto $\delta_A$. In applications of the coding lemma throughout the paper, the prewellordering and universal set used will always be the ones produced above. 

Therefore there is a $\Sigma_1$ formula $\varsigma(\alpha,e,A,\dot\reals)$ so that for all $\alpha < \delta_A$, $L(A,\dot\reals) \models \alpha \in S_e^{\rho_A} \Leftrightarrow \varsigma(\alpha,e,A,\dot\reals)$.
\end{definition}

\Begin{definition}{partition property}
A function $f : \lambda \rightarrow \mathrm{ON}$ has uniform cofinality $\omega$ if and only if there is a $g : \lambda \times \omega \rightarrow \mathrm{ON}$ with the property that for all $\alpha < \lambda$ and $n \in \omega$, $g(\alpha,n) < g(\alpha, n + 1)$ and $f(\alpha) = \sup\{g(\alpha,n) : n \in \omega\}$. 

A function $f : \lambda \rightarrow \mathrm{ON}$ is of the correct type if and only if $f$ is strictly increasing, for all $\alpha < \lambda$, $f(\alpha) > \sup\{f(\beta) : \beta < \alpha\}$, and $f$ has uniform cofinality $\omega$. 

Let $\kappa$ be an ordinal. For ordinals $\lambda \leq \kappa$, let $\kappa \rightarrow (\kappa)^\lambda_2$ denote that for all $\Phi : [\kappa]^\lambda \rightarrow 2$, there is a club $C \subseteq \kappa$ and $i \in 2$ so that for all $f : \lambda \rightarrow C$ of the correct type, $\Phi(f) = i$. 

If $\kappa \rightarrow (\kappa)^\kappa_2$, then one says that $\kappa$ has the strong partition property. If for all $\eta < \kappa$, $\kappa \rightarrow (\kappa)^\eta_2$, then $\kappa$ is said to have the weak partition property. (Note that $\kappa \rightarrow (\kappa)^2_2$ implies that $\kappa$ is regular.)
\end{definition}

\Begin{fact}{stable has strong partition}
Assume $\ZF + \AD$. Let $A \subseteq \reals$. Then $\delta_A$ has the strong partition property in $L(A,\reals)$ and even in $V$.
\end{fact}

\begin{proof}
This is shown by following Martin's template for establishing partition properties. The reflection properties and the uniform coding lemma is used to produce a good coding system for functions $f : \delta_A \rightarrow \delta_A$. See \cite{Structural-Consequences-AD} for more details. See \cite{Extensions-Axiom-Determinacy} for the details of this specific result.
\end{proof}

\Begin{definition}{wadge notions}
The ordinal $\Theta$ is the supremum of the ordinals which are surjective images of $\reals$. 

For $A,B \in \mathscr{P}(\reals)$, $A \leq_w B$ denotes that $A$ is Wadge reducible to $B$. For each $r \in \reals$, let $\Xi_r$ denote the Wadge reduction coded by $r$. So $\Xi_r^{-1}[B]$ is the subset of $\reals$ reducible to $B$ via the Wadge reduction coded by $r$.

The Wadge lemma states that $\ZF + \AD$ implies that for all $A, B \in \mathscr{P}(\reals)$ either $A \leq_w B$ or $B \leq_w (\reals \setminus A)$.
\end{definition}

\Begin{fact}{strong partition cardinal cofinal theta}
(\cite{Axiom-Determinacy-Strong-Partition-Nonsingular}) Assume $\mathsf{ZF + AD}$. For all $\lambda < \Theta$, there exists some $\kappa$ with $\lambda < \kappa < \Theta$ so that $\kappa$ has the strong partition property. 
\end{fact}

\begin{proof}
This result follows from Fact \ref{stable has strong partition}. \cite{Axiom-Determinacy-Strong-Partition-Nonsingular} works with $\mathsf{ZF + DC + AD}$ as its base theory. \cite{Extensions-Axiom-Determinacy} has a careful presentation of this result from just $\ZF + \AD$.
\end{proof}

\Begin{definition}{ground club property}
Let $\kappa$ be a regular cardinal and $\bbP = (\bbP,\leq_\bbP,1_{\bbP})$ be a forcing. $\bbP$ has the ground club property at $\kappa$ if and only if for all $p \in \bbP$ and all $\bbP$-name $\dot D$ such that $p \forces_{\bbP}$ ``$\dot D$ is a club subset of $\check \kappa$'', there is some club $C \subseteq \kappa$ so that $p \forces_\bbP \check C \subseteq \dot D$. 
\end{definition}

\Begin{lemma}{strong partition ground club same reals}
Assume $\ZF$. Let $\bbP$ be a forcing and $p \in \bbP$. If $\bbP$ has the ground club property at $\kappa$ and $p \forces \kappa \rightarrow (\kappa)^\omega_2$, then $p \forces \dot \reals = \check \reals$. 
\end{lemma}

\begin{proof}
Let $G \subseteq \bbP$ be any $\bbP$-generic filter over $V$ containing $p$. Observe that every set of $V$ belongs to $V[G]$, so in particular, $([\kappa]^\omega)^V \in V[G]$.

In $V[G]$, define $\Phi : [\kappa]^\omega \rightarrow 2$ by
$$\Phi(f) = \begin{cases}
0 & \quad f \in ([\kappa]^\omega)^V \\
1 & \quad \text{otherwise}
\end{cases}.$$

Let $D \subseteq \kappa$ be a club set homogeneous for $\Phi$. By the ground club property at $\kappa$, there is some $C \subseteq D$ with $C \in V$ and is a club in $V$. Pick any $f \in ([C]^\omega)^V$ of correct type. Then $\Phi(f) = 0$. Thus $D$ is homogeneous for $\Phi$ taking value $0$. Therefore $C$ is also homogeneous for $\Phi$ taking value $0$. Any function $f \in ([C]^\omega)^{V[G]}$ of the correct type belongs to $V$.

Let $c_i = C(\omega \cdot i + \omega)$. Since $C \in V$, the sequence $(c_i : i \in \omega)$ belongs to $V$. Each $c_i \in C$ since $C$ is club and each $c_i$ has cofinality $\omega$. Let $z \in \reals^{V[G]}$. Let $f_z = \{c_i : i \in z\}$. Then $f_z \in [C]^\omega$ and is of correct type. So $\Phi(f_z) = 0$. $f_z \in V$. Then $z = \{i \in \omega : c_i \in f_z\}$. So $z \in \reals^V$. 
\end{proof}

\section{Wellorderable Forcings of Cardinality Less than $\Theta$}\label{wo forcing less than theta}

This section will show that a nontrivial forcing on a wellorderable set of cardinality less than $\Theta$ can not preserve $\AD$. The results of this section are subsumed by the results of Section \ref{destroying AD Theta regular}; however, the argument there is far less natural for wellorderable forcings.

\Begin{fact}{ground club property above size of forcing}
Assume $\ZF$. Let $\bbP$ be a wellorderable forcing of size $\lambda$. Then $\bbP$ has the ground club property at $\kappa$ for all regular $\kappa > \lambda$.
\end{fact}

\begin{proof}
Let $p \in \bbP$ and $\dot D$ be a $\bbP$ name such that $p \forces_\bbP$ ``$\dot D \subseteq \check \kappa$ is a club''. For each $\alpha < \kappa$, let $A_\alpha = \{p \in \bbP : (\exists \beta < \kappa)(p \forces_\bbP \dot D(\check \alpha) = \check \beta)\}$. Let $B_\alpha = \{\beta : (\exists p \in A_\alpha)(p \forces_{\bbP} \dot D(\check \alpha) = \check \beta)\}$. Since $|A_\alpha| \leq |\bbP| = \lambda < \kappa$ and $\kappa$ is regular, $\sup B_\alpha < \kappa$. Let $F(\alpha) = \sup B_\alpha$. Note that $F(\alpha) \geq \alpha$ since $1_{\bbP} \forces_{\bbP} \dot D(\check \alpha) \geq \check \alpha$. Let $C = \{\alpha < \kappa : (\forall \eta < \alpha)(F(\eta) < \alpha)\}$. $C$ is club.

Let $G \subseteq \bbP$ be a $\bbP$-generic filter over $V$ with $p \in G$. Let $D = \dot D[G]$. Suppose $\alpha \in C$. Since $G$ is generic, for each $\eta < \alpha$, $G \cap A_\eta \neq \emptyset$. For any $q \in G \cap A_\eta$, $q \forces_{\bbP} \dot D(\check \eta) < F(\check \eta) < \check \alpha$. Hence $\eta \leq D(\eta) < \alpha$ for all $\eta < \alpha$. Since $D$ is a club, $\alpha \in D$. This shows $C \subseteq D$ in $V[G]$. Since $G$ was arbitrary with $p \in G$, $p \forces_\bbP \check C \subseteq \dot D$. 
\end{proof}

\Begin{theorem}{WO forcing below theta destroy AD}
Assume $\mathsf{ZF + AD}$. If $\bbP$ is a nontrivial wellorderable forcing of cardinality less than $\Theta$, then $1_{\bbP} \forces_{\bbP} \neg \AD$. 
\end{theorem}

\begin{proof}
Suppose $|\bbP| = \delta$ where $\delta < \Theta$ is a cardinal. One may assume $\bbP \subseteq \delta$. 

Let $G \subseteq \bbP$ be a $\bbP$-generic filter over $V$. Assume that $V[G] \models \AD$. By Fact \ref{strong partition cardinal cofinal theta}, let $\kappa$ be a cardinal such that $\delta < \kappa < \Theta^{V[G]}$ and has the strong partition property in $V[G]$. Therefore, $\kappa$ is regular in $V$. Let $p \in G$ be such that $p \forces \kappa \rightarrow (\kappa)^\omega_2$. Since $\delta < \Theta$, let $\pi : \reals^{V} \rightarrow \delta$ be a surjection in $V$.

Since $\bbP \subseteq \delta$, if $\bbP$ is nontrivial, then $G$ is a new subset of $\delta$. Since $V[G] \models \AD$, there is some $e \in \reals^{V[G]}$ so that $S_e^{\pi} = G$ by Fact \ref{coding sets of ordinals}. If $\reals^V = \reals^{V[G]}$, then this would imply $G \in V$. Hence one must have that $\reals^V \subsetneq \reals^{V[G]}$. 

Fact \ref{ground club property above size of forcing} implies that $\bbP$ has the ground club propery at $\kappa$. Lemma \ref{strong partition ground club same reals} implies that $p \forces_\bbP \dot \reals = \check \reals$. So $\reals^V = \reals^{V[G]}$. Contradiction.
\end{proof}

The previous theorem illustrates the main ideas to be used in Section \ref{destroying AD Theta regular}. The above proof uses the partition property $\kappa \rightarrow (\kappa)^\omega_2$. This requires the theorem to be restricted to wellorderable forcings of cardinality less than $\Theta$. In discussion with Goldberg, the following more elementary argument was found which could apply to more wellorderable forcings:

\Begin{fact}{wellorderable forcing adding new reals destroy perfect set prop}
$(\ZF)$ Assume all sets of reals have the Baire property. Let $\bbP$ be a wellorderable forcing such that $1_\bbP \forces_\bbP \check \reals \subsetneq \dot \reals$ (adds new reals), then $1_\bbP \forces_\bbP$ ``$\check \reals$ has no perfect subset''.
\end{fact}

\begin{proof}
Suppose there was a $G \subseteq \bbP$ which is $\bbP$-generic over $V$ and $V[G] \models \reals^V$ has a perfect subset. In $V[G]$, let $T$ be a perfect tree so that $[T] \subseteq \reals^V$. Let $\dot T$ be a name for $T$ and $q \in G$ be such that $q \forces_\bbP \dot T$ is a perfect tree.

Work in $V$. For each $p \in \bbP$, let $A_p = \{x \in \reals : p \forces_\bbP \check x \in [\dot T]\}$. Note that if $p \leq_\bbP q$, then each $A_p$ is closed. To see this: Suppose $z$ is a limit point of $A_p$. Let $H$ be any $\bbP$-generic filter over $V$ containing $p$. Since $[\dot T[H]]$ is a closed set and $A_p \subseteq [\dot T[H]]$, $z \in [\dot T[H]]$. Since $H$ was arbitrary containing $p$, $p \forces_\bbP \check z \in [\dot T]$. 

Note that in $V[G]$, $[T] \subseteq \bigcup_{p \leq_\bbP q} A_p$. Thus in $V$, $\bigcup_{p \leq_\bbP q} A_p$ is an uncountable set. By the Baire property in $V$ for all sets of reals, a wellordered union of meager sets is meager. Hence, there is some $p \in \bbP$ so that $A_p$ is uncountable. Since $A_p$ is a closed uncountable set, there is some perfect tree $U$ so that $[U] \subseteq A_p$. Note that for all $t \in U$, $p \forces_\bbP \check t \in \dot T$. Thus $p \forces_\bbP [\check U] \subseteq [\dot T]$. 

In $V[G]$, since $U \in V$ and $V[G]$ has a new real, $[U]$ must have a new real. Then $[\dot T[G]] = [T]$ has a new real. But $[T] \subseteq \reals^V$. Contradiction.
\end{proof}

\Begin{fact}{WO forcing can not WO R}
$(\ZF)$. Let $\bbP$ be a forcing on a wellorderable set. If $\reals$ is not wellorderable, then $1_\bbP \forces_\bbP \check\reals$ is not wellorderable.
\end{fact}

\begin{proof}
Since $\bbP$ is wellorderable, let $|\bbP| = \delta$ where $\delta$ is some ordinal. One may assume $\bbP \subseteq \delta$. Suppose $G \subseteq \bbP$ is $\bbP$-generic over $V$ and $V[G] \models \reals^V$ is wellorderable. There is an injection $\Phi : \reals^V \rightarrow \mathrm{ON}$. Let $\dot \Phi$ be a $\bbP$-name for $\Phi$. 

Work in $V$: For each $r \in \reals$, let $A_r = \{\langle p,\beta\rangle : p \forces_\bbP \dot \Phi(\check r) = \check \beta\}$, where $\langle \cdot, \cdot \rangle$ denotes a definable bijection of $\ON \times \ON$ with $\ON$. Each $A_r \neq \emptyset$ and if $r \neq s$, then $A_r \cap A_s = \emptyset$. In $V$, let $\Psi : \reals \rightarrow \mathrm{ON}$ be defined by $\Psi(r) = \min A_r$. $\Psi$ is an injection and hence $\reals^V$ is wellorderable in $V$. Contradiction.
\end{proof}

\Begin{corollary}{WO forcing adding new real destroy AD}
Assume $\ZF + \AD$. If $\bbP$ is a wellorderable forcing which adds a new real, then $1_\bbP \forces_\bbP \neg \AD$. 
\end{corollary}

\begin{proof}
Let $G \subseteq \bbP$ be $\bbP$-generic over $V$. By Fact \ref{WO forcing can not WO R}, $V[G]$ must think that $\reals^V$ is uncountable. By Fact \ref{wellorderable forcing adding new reals destroy perfect set prop}, $\bbR^V$ is an uncountable set of reals without the perfect set property. Thus $\AD$ must fail.
\end{proof}

\Begin{question}{WO forcing adds new reals}
Assume $\ZF + \AD$. Can a nontrivial wellorderable forcing preserve $\AD$? 

If $\bbP$ is a nontrivial wellorderable forcing, then must $\bbP$ add a new real?

The proof of Theorem \ref{WO forcing below theta destroy AD} used the Moschovakis coding lemma to show that nontrivial wellorderable forcing of cardinality less than $\Theta$ must add a new real.
\end{question}

\section{Preservation of $\Theta$}\label{preservation Theta}

Trang and Ikegami showed that in natural models of $\AD^+$, every forcing that preserves $\AD$ must preserve $\Theta$:

\Begin{fact}{natural model AD+ preserve theta}
(Ikegami and Trang) Assume $\mathsf{ZF + AD^+ + V = L(\mathscr{P}(\reals))}$. If $\bbP$ is a nontrivial forcing and $1_{\bbP} \forces_{\bbP} \AD$, then $1_{\bbP} \forces \dot\Theta = \Theta^V$.
\end{fact}

This section will show under $\ZF + \AD$ that any forcing which is a surjective image of $\reals$ that preserves $\AD$ must preserve $\Theta$. It will first be shown using Lemma \ref{strong partition ground club same reals} that any forcing that adds a new real and preserves $\AD$ must preserve $\Theta$. 

A nontrivial forcing adds the generic filter as a new object. If $\bbP$ is a surjective image of $\reals$, then a new set of reals must be added. It will then be shown under $\ZF + \AD$ that any nontrivial forcing which is a surjective image of $\reals$ which preserves $\AD$ must actually add a new real. Hence any nontrival forcing which is a surjective image of $\reals$ must preserve $\Theta$. 

Lemma \ref{real sized forcing preserve Theta} and Fact \ref{nontrivial forcing adds new reals} below have been known to Ikegami and Trang under $\mathsf{ZF + AD^+ + V = L(\mathscr{P}(\reals))}$ for forcing more general than those which are surjective images of $\reals$. An important aspect of their argument involves the sharps of sets of reals. It should be noted that the arguments below are for forcing which are surjective images of $\reals$ proved under just $\ZF + \AD$ without $\DC_\reals$. $\DC_\reals$ is used in some classical arguments to produce sharps of sets of reals and to show the wellfoundedness of the Wadge hierarchy.

\Begin{fact}{R sized forcing ground club above Theta}
Let $\bbP$ be a forcing which is a surjective image of $\reals$. For each regular $\kappa \geq \Theta$, $\bbP$ has the ground club property at $\kappa$. 
\end{fact}

\begin{proof}
Let $\pi : \reals \rightarrow \bbP$ be a surjection. Let $\kappa \geq \Theta$ be regular.

Let $p \in \bbP$ and $\dot D$ be a $\bbP$-name so that $p \forces_\bbP$ ``$\dot D \subseteq \check \kappa$ is a club''. For each $\alpha < \kappa$, let $A_\alpha = \{p \in \bbP : (\exists \beta < \kappa)(p \forces_\bbP \dot D(\check \alpha) = \check \beta)\}$. Let $B_\alpha = \{\beta : (\exists p \in A_\alpha)(p \forces \dot D(\check \alpha) = \check \beta)\}$. Define in $V$, $\Phi: \reals \rightarrow \kappa$ by
$$\Phi(r) = \begin{cases}
0 & \quad \pi(r) \notin A_\alpha \\
\beta & \quad \pi(r) \in A_\alpha \wedge \pi(r) \forces_\bbP \dot D(\check \alpha) = \check \beta 
\end{cases}$$
$\Phi$ induces a prewellordering on $\reals$. Let $\delta < \Theta^V$ be the length of this prewellordering. Hence $\Phi$ induces a map $\Psi : \delta \rightarrow \kappa$. Since $\kappa$ is regular in $V$, $\Psi$ must be bounded below $\kappa$.

Thus $\sup B_\alpha < \kappa$. Let $F(\alpha) = \sup B_\alpha$. Let $C = \{\alpha < \kappa : (\forall \eta < \alpha)(F(\eta) < \alpha)\}$. $C$ is a club subset of $\kappa$ in $V$. As in the proof of Fact \ref{ground club property above size of forcing}, $p \forces_\bbP \check C \subseteq \dot D$.
\end{proof}

\Begin{lemma}{real sized forcing preserve Theta}
Assume $\mathsf{ZF + AD}$. If $\bbP$ is a forcing which is a surjective image of $\reals$ and adds a new real, then $1_\bbP \forces_\bbP \AD$ implies that $1_{\bbP} \forces_\bbP \Theta = \Theta^V$. 
\end{lemma}

\begin{proof}
Let $G \subseteq \bbP$ be a $\bbP$-generic filter over $\bbP$. Suppose $V[G] \models \AD$ and $\Theta^{V[G]} > \Theta^V$. 

By Fact \ref{strong partition cardinal cofinal theta} applied in $V[G]$, there is a $\kappa$ such that $\Theta^V < \kappa < \Theta^{V[G]}$ and $\kappa \rightarrow (\kappa)^\omega_2$. Note $\kappa \rightarrow (\kappa)^2_2$ implies that $\kappa$ is regular in $V[G]$. Hence $\kappa$ is regular in $V$. By Fact \ref{R sized forcing ground club above Theta}, $\bbP$ has the ground club property at $\kappa$. Choose $p \in G$ so that $p \forces_\bbP \kappa \rightarrow (\kappa)^\omega_2$. Lemma \ref{strong partition ground club same reals} implies that $\reals^V = \reals^{V[G]}$. Contradiction.
\end{proof}

\Begin{fact}{nontrivial forcing adds new reals}
Assume $\ZF + \AD$. Let $\bbP$ be a nontrival forcing which is a surjective image of $\reals$. Suppose $1_\bbP \forces_\bbP \AD$. Then $1_{\bbP} \forces_\bbP \check\reals \subsetneq \dot\reals$. Hence $1_\bbP \forces_\bbP \dot\Theta = \Theta^V$.
\end{fact}

\begin{proof}
Let $\pi : \reals \rightarrow \bbP$ be a surjection. Suppose there is some $p \in \bbP$ so that $p \forces_\bbP \check \reals = \dot \reals$. Since $\bbP$ is a nontrivial forcing, $\pi^{-1}[\dot G]$ is forced to be a new set of reals. Since $p \forces_\bbP \check \reals = \dot \reals$, for each $A \in \mathscr{P}(\reals)^V$, $p \forces_\bbP \check A \leq_w \pi^{-1}[\dot G]$.

In $V$, define $\Phi : \reals \times \reals \rightarrow \Theta$ by
$$\Phi(r,s) = \begin{cases}
\alpha & \quad \text{$\pi(r) \forces_\bbP ``\Xi_s^{-1}[\pi^{-1}[\dot G]] \in \check V$ and is a prewellordering of length $\check \alpha$''} \\
0 & \quad \text{otherwise}
\end{cases}$$
Thus in $V$, $\Phi$ is a surjection of $\reals \times \reals$ onto $\Theta$. This is impossible.
\end{proof}

\Begin{fact}{surjective image R and forcing on R}
Assume that $\bbP$ is a forcing which is a surjective image of $\reals$. Then there is a forcing $\bbQ$ on $\reals$ so that for every $G \subseteq \bbP$ which is $\bbP$-generic over $V$, there is an $H \subseteq \bbQ$ which is $\bbQ$-generic over $V$ so that $V[G] = V[H]$. 
\end{fact}

\begin{proof}
Let $\pi : \reals \rightarrow \bbP$ be a surjection. Define a forcing $\bbQ$ on $\reals$ by $p \leq_\bbQ q$ if and only if $\pi(p) \leq_\bbP \pi(q)$. If $G \subseteq \bbP$ is a $\bbP$-generic filter over $V$, then $\pi^{-1}[G] \subseteq \bbQ$ is a $\bbQ$-generic filter over $V$ and $V[G] = V[\pi^{-1}[G]]$. 
\end{proof}

\Begin{lemma}{preserving AD and preserving V = L(R)}
Assume $\ZF + \AD$ and there is an $A \subseteq \reals$ such that $V = L(A,\reals)$. Let $\bbP$ be a forcing on $\reals$ such that $1_\bbP \forces_\bbP \AD$ and $\bbP \leq_w A$. Let $A \oplus \reals^V$ indicate some fixed recursive coding of the two sets of reals into a single set of reals. (Note that $V = L(A \oplus \reals^V,\reals^V)$.) Then $1_\bbP \forces \dot V = L(\check A\oplus \check \reals, \dot\reals)$. 
\end{lemma}

\begin{proof}
Suppose not. Let $G \subseteq \bbP$ be a $\bbP$-generic filter over $L(A,\reals)$ witnessing the failure of the conclusion of the lemma. Here $\reals$ refers to $\reals^{L(A,\reals)}$. Let $\reals^* = \dot\reals^{L(A,\reals)[G]}$. Note that $\bbP, \reals \in L(A \oplus \reals,\reals^*)$. Therefore, $L(A,\reals)$ is a definable inner model of $L(A \oplus\reals,\reals^*)$. Thus $\Theta^{L(A,\reals)} \leq \Theta^{L(A\oplus\reals,\reals^*)}$. Since $\bbP,\bbR\in L(A\oplus\reals,\reals^*)$, $L(A\oplus\reals,\reals^*) \neq L(A,\reals)[G]$ implies that $G \notin L(A\oplus\reals,\reals^*)$. Since $L(A\oplus\reals,\reals^*)$ and $L(A,\reals)[G]$ have the same set of reals, $G$ Wadge reduces every set of reals in $L(A,\reals^*)$.  In $L(A,\reals)[G]$, define $\Phi : \reals^* \rightarrow \Theta^{L(A,\reals^*)}$ by
$$\Phi(r) = \begin{cases}
\text{length}(\Xi_r^{-1}[G]) & \quad \text{$\Xi_r^{-1}[G] \in L(A\oplus\reals,\reals^*)$ and is a prewellordering on $\reals^*$} \\
0 & \quad \text{otherwise}
\end{cases}$$

$\Phi$ is a surjection in $L(A,\reals)[G]$ of $\reals^*$ onto $\Theta^{L(A \oplus \reals,\reals^*)}$. This implies that $\Theta^{L(A,\reals)} \leq \Theta^{L(A\oplus\reals,\reals^*)} < \Theta^{L(A,\reals)[G]}$. This contradicts Fact \ref{nontrivial forcing adds new reals} which asserts that $\Theta^{L(A,\reals)} = \Theta^{L(A,\reals)[G]}$. 
\end{proof}

\Begin{fact}{preserving AD implies preserve Theta is regular}
Assume $\ZF + \AD$. If $\bbP$ is a forcing which is the surjective image of $\reals$ and $\Theta$ is regular, then $1_\bbP\forces_{\bbP}\AD$ implies $1_\bbP \forces \Theta$ is regular.
\end{fact}

\begin{proof}
Let $\pi : \reals \rightarrow \bbP$ be a surjection. Let $G \subseteq \bbP$ be a $\bbP$-generic filter over $V$. By Fact \ref{nontrivial forcing adds new reals}, $V[G] \models \Theta^{V[G]} = \Theta^V$. Suppose $\Theta$ is not regular in $V[G]$. There is some $\eta < \Theta$ and a function $f : \eta \rightarrow \Theta$ which is cofinal. Let $\tau \in V$ be a $\bbP$-name so that $\tau[G] = f$. 

Now work in $V$. Define $g : \eta \times \reals \rightarrow \Theta$ by
$$g(\alpha,r) = \begin{cases}
0 & \quad (\forall \beta < \Theta)(\pi(r) \not\forces_\bbP \tau(\check \alpha) = \check\beta) \\
\beta & \quad \pi(r) \forces_\bbP \tau(\check \alpha) = \check \beta
\end{cases}.$$
Let $\rho : \reals \rightarrow \eta$ be a surjection. Define $h : \reals \rightarrow \Theta$ by $h(x) = g(\rho(x_1),x_2)$, where $x = \langle x_1,x_2\rangle$ under some standard pairing function. Let $x \preceq y$ if and only if $h(x) \leq h(y)$. As $\preceq$ is a prewellordering of $\reals$, it has length some $\delta < \Theta$. Thus there is a map $\tilde h : \delta \rightarrow \Theta$ which is cofinal. This is impossible since $\Theta$ is regular in $V$.
\end{proof}

\section{Destroying $\AD$ When $\Theta$ Is Regular}\label{destroying AD Theta regular}

By Fact \ref{surjective image R and forcing on R}, this section will assume that the forcing is on $\reals$. For such a forcing $\bbP$, a name for a real consisting of elements of the form $(\check n, p)$ for $n \in \omega$ and $p \in \bbP$ can be considered subsets of $\reals$. In this section, when one writes that a name $\sigma \in \mathscr{P}(\reals)$, it is understood that $\sigma$ takes this form.

\Begin{definition}{name condition}
Let $\bbP$ be a forcing on $\reals$. $\bbP$ has the name condition if and only if there is an $A \subseteq \reals$ so that $\bbP \leq_w A$ and $1_\bbP \forces_\bbP$ ``for all $r \in \dot\reals$, there is a $\bbP$-name $\sigma \in \mathscr{P}(\check \reals)^{L(\check A,\check \reals)}$ so that $\sigma[\dot G] = r$ and $L(\check A, \check \reals) \models \sigma \leq_w \check A$''.

This means that there is a set $A \subseteq \reals$ so that for all $G \subseteq \bbP$ which are $\bbP$-generic over $V$, for all $r \in \reals^{V[G]}$, there is a set of reals $\sigma$ in $L(A,\reals)$ which is also Wadge reducible to $A$ in $L(A,\reals)$ so that when $\sigma$ is construed as a $\bbP$-name, $\sigma[G] = r$. 
\end{definition}

\Begin{fact}{wellorderable forcing size less than theta name condition}
Assume $\ZF + \AD$. Suppose $\bbP$ is a wellorderable forcing of cardinality less than $\Theta$. Then $\bbP$ has the name condition.
\end{fact}

\begin{proof}
Suppose $|\bbP| = \delta$ where $\delta < \Theta$. One may assume $\bbP \subseteq \delta$. Suppose $\tau$ is a $\bbP$ name so that for some $p \in \bbP$, $p \forces_\bbP \tau \in \dot \reals$. Let $\sigma = \{(\check n, q) : q \forces_\bbP \check n \in \tau\}$. Then $p \forces_\bbP \sigma = \tau$. Note that $\sigma$ can be identified as a subset of $\delta$. Since $\delta < \Theta$, let $\preceq$ be a prewellordering of rank $\delta$. By the Moschovakis coding lemma, every subset of $\delta$ is coded by a real using $\preceq$. Thus $\sigma \in L(\preceq,\reals)$.
\end{proof}

\Begin{fact}{pwo that a set can wadge reduce}
Assume $\ZF + \AD$. Let $A \subseteq \reals$. Let $C_A$ be the set of ordinals $\alpha$ less than $\Theta$ so that $A$ can Wadge reduce a prewellordering on $\reals$ of length $\alpha$. Then $C_A$ is bounded below $\Theta$.
\end{fact}

\begin{proof}
Suppose not. Define $\Psi : \reals \times \reals \rightarrow \Theta$ by 
$$\Psi(r,s) = \begin{cases}
\mathrm{rk}_{\Xi_r^{-1}[A]}(s) & \quad \text{if $\Xi_r^{-1}[A]$ is a prewellordering on $\reals$} \\
0 & \quad \text{otherwise}
\end{cases}$$
where if $\preceq$ is a prewellordering on $\reals$, then $\mathrm{rk}_\preceq(s)$ denote the rank of $s$ in the prewellordering $\preceq$. $\Psi$ is a surjection of $\reals \times \reals$ onto $\Theta$. Contradiction.
\end{proof}

\Begin{fact}{Theta regular name condition}
Assume $\ZF + \AD + \Theta$ is regular. Every forcing $\bbP$ on $\reals$ such that $1_\bbP \forces_\bbP \AD$ has the name condition.
\end{fact}

\begin{proof}
Let $p \in \bbP$ and $G \subseteq \bbP$ be a $\bbP$-generic filter over $V$ such that $p \in G$. By Fact \ref{nontrivial forcing adds new reals} and Fact \ref{preserving AD implies preserve Theta is regular}, $\Theta^{V[G]} = \Theta^{V}$ and $\Theta$ remains regular in $V[G]$. 

Suppose $r \in \reals^{V[G]}$. There is some $\bbP$-name $\tau \in V$ so that $r = \tau[G]$. Let $\sigma = \{(\check n, s) : s \forces_\bbP \check n \in \tau\}$. Note that $\sigma[G] = \tau[G]$ and $\sigma$ can be considered as essentially a set of reals.

Since $\mathscr{P}(\reals)^V \in V[G]$, one can define a function $\Phi : \reals^{V[G]} \rightarrow \Theta$ by 
$$\Phi(r) = \min\{\sup(C_\sigma)^V + 1 : \sigma \in \mathscr{P}(\reals)^V \wedge \sigma[G] =r \}$$
where $C_A$, for $A \subseteq \reals$, is defined in Fact \ref{pwo that a set can wadge reduce}.

In $V[G]$, define $x \sqsubseteq y$ if and only if $\Phi(x) \leq \Phi(y)$. $\sqsubseteq$ is a prewellordering on $\reals$. There is some $\delta < \Theta^{V[G]} = \Theta^V$ so that $\sqsubseteq$ has length $\delta$. Thus $\Phi$ induces a map $\tilde\Phi : \delta \rightarrow \Theta$. Since $\Theta$ is regular in $V[G]$, $\tilde \Phi$ and hence $\Phi$ is bounded below some $\gamma < \Theta$. 

Fix a prewellordering $\preceq^*$ in $V$ of length greater than or equal to $\gamma$. Let $r \in \reals^{V[G]}$. Let $\sigma \in V$ be a set of reals so that when it is construed as a $\bbP$-name, $\sigma[G] = r$ and $\Phi(r) = (C_\sigma)^V + 1$. Since $\gamma > \sup (C_\sigma)^V$, $\sigma$ can not Wadge reduce $\preceq^*$ in $V$. Hence by Wadge's lemma, $\sigma \leq_w \preceq^*$ in $V$. 

It has been shown that in $V[G]$, there is some ordinal $\gamma$, so that for any prewellordering $\preceq^* \in V$ of length greater than or equal to $\gamma$, every $r \in \reals^{V[G]}$ has a name $\sigma \in \mathscr{P}(\reals)^{L(\preceq^*,\reals)}$ so that $\sigma[G] = r$ and $L(\preceq^*,\reals) \models \sigma \leq_w \preceq^*$. Find some $q \leq_\bbP p$, $q \in G$, and some $\gamma < \Theta$ so that $q$ which forces this above statement about $\gamma$. Since $p \in \bbP$ was arbitrary, it has been shown that there is a dense set of $q$ for which there is some $\gamma$ so that $q$ forces the above statement involving $\gamma$.

Define $\Psi : \bbP \rightarrow \Theta$ by $\Psi(q)$ is the least $\gamma$ so that $q$ forces the above statement involving $\gamma$ if such a $\gamma$ exists. Let $\Psi(q) = 0$ otherwise. $\Psi$ induces a prewellordering on $\reals$ of length $\delta < \Theta$. Since $\Theta$ is regular in $V$, $\Psi$ is bounded below $\Theta$ by some $\gamma$. Let $\preceq^*$ be some prewellordering on $\reals$ of length $\gamma$. Let $A = \preceq^*$. One has that $A$ witnesses that $\bbP$ has the name condition.
\end{proof}

\Begin{lemma}{name condition implies ground club}
Assume $\ZF + \AD$. Let $\bbP$ be a forcing on $\reals$ and $1_\bbP \forces \AD$. Assume that $\bbP$ has the name condition. Let $A \subseteq \reals$ witness the name condition. Then in $L(A,\reals)$, $1_\bbP \forces_\bbP \AD$, $\delta_A$ has the ground club property, and $1_\bbP \forces_\bbP \check\delta_A$ has the strong partition property.
\end{lemma}

\begin{proof}
Let $A$ witness the name condition. Note that $L(A,\reals) \models \AD$. 

Throughout this proof, $\reals$ denotes $\reals^V$ and $\reals^*$ denotes $\reals^{V[G]}$ whenever $G$ is $\bbP$-generic over $V$.

Let $p \in \bbP$. Let $G \subseteq \bbP$ be any $\bbP$-generic filter over $V$ containing $p$. By definition of the name condition, $\dot\reals^{L(A,\reals)[G]} = \dot\reals^{V[G]}$. Thus since $V[G] \models \AD$, $L(A,\reals)[G] \models \AD$. Let $q \leq_\bbP p$ with $q \in G$ be such that $L(A,\reals) \models q \forces_\bbP \AD$. Since $p \in \bbP$ was arbitrary, there is a dense set of $q \in \bbP$ so that $L(A,\reals) \models q \forces_\bbP \AD$. One has that $L(A,\reals) \models 1_\bbP \forces_\bbP \AD$. 

By Lemma \ref{preserving AD and preserving V = L(R)}, $L(A,\reals)[G] = L(A \oplus \reals,\reals^*)$. 

Let $p \in G$. Let $G \subseteq \bbP$ be any $\bbP$-generic filter over $V$.

Claim 1: $\delta_A = (\delta_{A\oplus\reals})^{L(A\oplus\reals,\reals^*)}$.

Let $r \in \reals^*$. By the name condition, there is some $\tau \subseteq \reals$ which is Wadge reducible to $A$ and $\tau[G] = r$ when $\tau$ is construed as a $\bbP$-name. Note that every set which is Wadge reducible to $A$ appears at level $L_1(A,\reals)$. Let $\varphi(\dot v, A\oplus \reals, \dot \reals)$ be a $\Sigma_1$ formulas. Suppose that $L(A\oplus\reals,\reals^*) \models \varphi(r,A\oplus\reals,\reals^*)$. Since $L(A\oplus\reals,\reals^*) = L(A,\reals)[G]$, there is some $q_0 \leq_\bbP p$ so that $q_0 \in G$ and
$$L(A,\reals) \models q_0 \forces_\bbP L(\check A \oplus \check \reals,\dot\reals) \models \varphi(\tau,\check A\oplus \check \reals, \dot \reals).$$
By replacement, the following is a true $\Sigma_1(L(A,\reals), \reals \cup \{\reals, A\})$ formula: (Note that it is important that $\tau \leq_w A$.)
$$L(A,\reals) \models (\exists \alpha)(L_\alpha(A,\reals) \models q_0 \forces_\bbP L(\check A\oplus\check\reals,\dot\reals) \models \varphi(\tau,\check A\oplus\check\reals, \dot\reals)).$$
By definition of $\delta_A$, there exists some $\alpha < \delta_A$ so that
$$L(A,\reals) \models L_\alpha(A,\reals) \models q_0 \forces_\bbP L(\check A\oplus\check\reals,\dot\reals)\models \varphi(\tau,\check A\oplus\check\reals,\dot\reals).$$
Hence for some $\alpha < \delta_A$,
$$L_\alpha(A,\reals) \models q_0 \forces_\bbP L(\check A\oplus\check\reals,\dot\reals)\models \varphi(\tau,\check A\oplus\check\reals,\dot\reals).$$
Since $q_0 \in G$, the forcing theorem gives
$$L_\alpha(A,\reals)[G] \models L(A\oplus\reals,\dot \reals) \models \varphi(r,A\oplus\reals,\dot \reals).$$
Also 
$$(L(A\oplus\reals,\dot\reals))^{L_\alpha(A,\reals)[G]} = L_\alpha(A\oplus\reals,\dot\reals^{L_\alpha(A,\reals)[G]}).$$
Since $A$ witnesses the name condition, every $t \in \reals^*$ has a name in $L_1(A,\reals)$. Hence $\dot\reals^{L_\alpha(A,\reals)[G]} = \reals^*$. Thus one has
$$L_\alpha(A,\reals)[G] \models L_\alpha(A\oplus\reals,\reals^*) \models \varphi(r,A\oplus\reals,\reals^*).$$
Thus
$$L_\alpha(A\oplus\reals,\reals^*) \models \varphi(r,A\oplus\reals,\reals^*).$$
By upward absolute of $\Sigma_1$ formulas,
$$L_{\delta_A}(A\oplus\reals,\reals^*) \models \varphi(r,A\oplus\reals,\reals^*).$$
It has been established that $(\delta_{A\oplus\reals})^{L(A\oplus\reals,\reals^*)} \leq \delta_A$.

Let $\varphi(\dot v, A, \dot \reals)$ be a $\Sigma_1$ formula and $r \in \reals^V$. Note $\bbR \in L(A\oplus\reals,\reals^*)$. Suppose $L(A,\reals) \models \varphi(r,A,\reals)$. Then 
$$L(A\oplus\reals,\reals^*) \models L(A,\reals) \models \varphi(r,A,\reals).$$ 
The following is a true $\Sigma_1(L(A\oplus\reals,\reals^*),\reals^*\cup\{A\oplus\reals,\reals^*\})$ sentence
$$L(A\oplus\reals,\reals^*) \models (\exists \alpha)(L_\alpha(A,\reals) \models \varphi(r,A,\reals)).$$
By definition of $(\delta_{A\oplus\reals})^{L(A\oplus\reals,\reals^*)}$, there is some $\alpha < (\delta_{A\oplus\reals})^{L(A\oplus\reals,\reals^*)}$ so that
$$L(A\oplus\reals,\reals^*) \models L_\alpha(A,\reals) \models \varphi(r,A,\reals).$$
Thus
$$L_\alpha(A,\reals) \models \varphi(r,A,\reals).$$
By upward absoluteness
$$L_{(\delta_{A\oplus\reals})^{L(A\oplus\reals,\reals^*)}}(A,\reals) \models \varphi(r,A,\reals).$$
This shows that $\delta_A \leq (\delta_{A\oplus\reals})^{L(A\oplus\reals,\reals^*)}$. Claim 1 has been established.

By Claim 1, let $q \leq_\bbP p$ with $q \in G$ be such that $L(A,\reals) \models q \forces \check\delta_A = \dot\delta_{A\oplus\check\reals}$. Since $p \in \bbP$ was arbitrary, the set of $q \in \bbP$ such that $L(A,\reals) \models q \forces \check \delta_A = \dot \delta_{A \oplus \check \reals}$ is dense. Thus $L(A,\reals) \models 1_\bbP \forces_\bbP \check \delta_A = \dot \delta_{A\oplus\check\reals}$. 

Fact \ref{stable has strong partition} now gives that $1_\bbP \forces_\bbP \check\delta_A$ has the strong partition property. It remains to show that $\delta_A$ has the ground club property.

Claim 2: In $L(A,\reals)$, $\delta_A$ has the ground club property.

Let $p \in \bbP$ and $G$ be $\bbP$-generic over $L(A,\reals)$ containing $p$. Recall again that $L(A,\reals)[G] = L(A \oplus \reals, \reals^*)$ by Lemma \ref{preserving AD and preserving V = L(R)} and $\delta_A^{L(A,\reals)} = \delta_{A \oplus \reals}^{L(A \oplus \reals, \reals^*)}$. Let $D \in L(A,\reals)[G]$ be a club subset of $\delta_A^{L(A,\reals)} = \delta_{A\oplus\reals}^{L(A\oplus\reals,\reals^*)}$. Let $\rho_{A \oplus \reals}$ and $\varsigma$ be those objects from Definition \ref{prewellordering of stable} for $A \oplus \reals$ defined in $L(A,\reals)[G] = L(A \oplus \reals,\reals^*)$. Since $L(A,\reals)[G] = L(A\oplus\reals,\reals^*) \models \AD$ and Fact \ref{coding sets of ordinals}, there is some $e \in \reals^*$ so that the graph of the increasing enumeration of $D$ is $S_e^{\rho_{A\oplus\reals}}$. By the name condition as witnessed by $A$, there is some $\bbP$-name $\dot e \subseteq \reals$ so that $\dot e \leq_w A$ by a Wadge reduction coded in $L(A,\reals)$ and $\dot e[G] = e$. There is some $q_0 \leq_\bbP q$ with $q_0 \in G$ so that $q_0 \forces ``S_{\dot e}^{\rho_{A\oplus\reals}}$ is the graph of an enumeration of a club subset of $\check\delta_A$''.

By reflection, for each $\beta < \delta_A$, the following is a true $\Sigma_1$ statement in $L(A,\reals)$ using parameters among $A$, $\reals$, and elements of $L_{\delta_A}(A,\reals)$:
$$L(A,\reals) \models (\exists \alpha)(L_\alpha(A,\reals) \models (\forall k \leq_\bbP q_0)(\exists j \leq_\bbP k)(\exists \gamma)(j \forces_\bbP \varsigma(\langle \check \beta,\check \gamma\rangle, \dot e, \check A \oplus \check \reals,\dot\reals)),$$
where $\langle \cdot,\cdot\rangle$ refers to a fixed ordinal pairing function. This merely states that there is a dense set of conditions below $q_0$ which forces a value for the image of $\check\beta$ under the function whose graph is $S^{\rho_{\check A\oplus\check\reals}}_{\dot e}$.

By the definition of $\delta_A$ in $L(A,\reals)$, there is some $\alpha < \delta_A$ so that
$$L_\alpha(A,\reals) \models (\forall k \leq_\bbP q_0)(\exists j \leq_\bbP k)(\exists \gamma)(j \forces_\bbP \varsigma(\langle \check \beta,\check \gamma\rangle, \dot e, \check A \oplus \check \reals,\dot\reals).$$
Let $\epsilon_\beta$ be the least $\alpha$ with this property. By upward absoluteness of the $\Sigma_1$ formula $\varsigma$,
$$L(A,\reals) \models q_0 \forces_\bbP (\exists \gamma < \check \epsilon_\beta)(\varsigma(\langle \check \beta,\gamma\rangle,\dot e,\check A \oplus \check \reals,\dot \reals)).$$
Thus for all $\bbP$-generic filter $H$ containing $q_0$, the $\beta^\text{th}$ element of the club subset of $\delta_A^{L(A,\reals)}$ enumerated by the function whose graph is $S_{\dot e[H]}^{\rho_{A\oplus\reals}}$ is less than $\epsilon_\beta$. Define in $L(A,\reals)$, a function $g : \delta_A \rightarrow \delta_A$ by $g(\beta) = \epsilon_\beta$. Let $C = \{\mu < \delta_A : (\forall \gamma < \mu)(g(\gamma) < \mu)\}$. By the same argument as in the proof of Fact \ref{ground club property above size of forcing}, $C \subseteq \delta_A$ is a club in $L(A,\reals)$ and $q_0 \forces ``\check C$ is a subset of the club enumerated by $S_{\dot e}^{\rho_{A\oplus\reals}}$''. Thus $L(A,\reals)[G] \models C \subseteq D$. This proves Claim 2 and completes the lemma.
\end{proof}

\Begin{theorem}{destruction AD when theta regular}
Assume $\ZF + \AD +$ $\Theta$ is regular. Suppose $\bbP$ is a nontrivial forcing which is a surjective image of $\reals$. Then $1_\bbP \forces_\bbP \neg\AD$. 
\end{theorem}

\begin{proof}
By Fact \ref{surjective image R and forcing on R}, one may assume $\bbP \subseteq \reals$. Assume $\AD$ is preserved by the forcing. Fact \ref{Theta regular name condition} implies that $\bbP$ has the name condition. Let $A \subseteq \reals$ witness the name condition.

Work in $L(A,\reals)$. Fact \ref{nontrivial forcing adds new reals} states that a new real must be added. However Lemma \ref{strong partition ground club same reals} and Lemma \ref{name condition implies ground club} imply that the ground model and the forcing extension have the same reals. Contradiction.
\end{proof}

\Begin{corollary}{destruction AD natural model AD+}
Assume $\mathsf{ZF + AD + V = L(\reals)}$. No nontrivial forcing $\bbP \in L_\Theta(\reals)$ can preserve $\AD$.

In fact, assume $\mathsf{ZF + AD^+ + \neg\AD_\reals + V = L(\mathscr{P}(\reals))}$. No nontrivial forcing which is the surjective image of $\reals$ can preserve $\AD$.
\end{corollary}

\begin{proof}
If there is some set $X$ so that every set is $\OD_{X,r}$ for some $r \in \reals$, then $\Theta$ is regular. Hence if $L(\reals) \models \AD$, then $L(\reals) \models \Theta$ is regular. Woodin showed that if $\mathsf{ZF + AD^+ + \neg AD_\reals + V = L(\mathscr{P}(\reals))}$ holds, then there is some set of ordinals $J$ so that $V = L(J,\reals)$. Hence in these natural models of $\AD^+ + \neg\AD_\reals$, $\Theta$ is regular. 
\end{proof}

\Begin{question}{preservation of theta forcing on R}
Assume $\ZF + \AD$. If $\bbP$ is a nontrivial forcing which is a surjective image of $\reals$, then does $1_\bbP \forces_\bbP \neg \AD$ hold?

By the above, it remains to consider the case when $\Theta$ is singular.
\end{question}

Let $\Theta_0$ be the supremum of the ordinals which are the surjective image of $\reals$ by $\OD$ surjections. Ikegami and Trang have informed the authors that the consistency of $\ZF + \AD^+$ and $\Theta > \Theta_0$ implies the consistency of the statement that there is a forcing $\bbP$ (which is not a surjective image of $\reals$) such that $1_\bbP \forces_\bbP \AD \wedge \check \Theta < \dot \Theta$. This model also does not satisfy $\mathsf{ZF + \AD^+ + V = L(\mathscr{P}(\reals))}$.

\bibliographystyle{amsplain}
\bibliography{references}

\end{document}